\documentclass[11pt]{amsart} 
\usepackage{amsthm,amsbsy,amsmath,amssymb,amscd,amsfonts,array,mathrsfs,verbatim,enumerate,xypic,enumitem,color}
\usepackage[all]{xy}
\xyoption{arc} 

\newtheorem{thm}{Theorem}
     
\newtheorem{lem}[thm]{Lemma}
\newtheorem{cor}[thm]{Corollary}

\theoremstyle{definition}

\newtheorem{defn}{Definition}
\newtheorem{example}{Example} 
\newtheorem{rem}{Remark}

\DeclareFontFamily{OT1}{rsfs}{}
\DeclareFontShape{OT1}{rsfs}{n}{it}{<-> rsfs10}{}
\DeclareMathAlphabet{\curly}{OT1}{rsfs}{n}{it}

\newcommand{\SeqTop}{{\bf SeqTop}} 
\newcommand{\Sets}{{\bf Sets}} 
\newcommand{\Sch}{{\bf Sch}} 
\newcommand{\Top}{{\bf Top}}   

\newcommand{\defeq}{\mathrel{\vcenter{\baselineskip0.5ex \lineskiplimit0pt \hbox{\scriptsize.}\hbox{\scriptsize.}}} =}

\newcommand{\Z}{{\bf Z}} 

\newcommand{\NN}{\mathbb{N}} 
\newcommand{\RR}{\mathbb{R}} 
\newcommand{\QQ}{\mathbb{Q}} 
\newcommand{\PP}{\mathbb{P}}
\newcommand{\CC}{\mathbb{C}}

\newcommand{\UU}{\mathbb{U}}

\newcommand{\F}{\curly{F}} 
\newcommand{\M}{\mathcal{M}} 

\newcommand{\ov}{\overline}
\newcommand{\into}{\hookrightarrow}
\newcommand{\be}{\begin{eqnarray*}}
\newcommand{\ee}{\end{eqnarray*}}

\newcommand{\bne}[1]{\begin{eqnarray} \label{#1} }
\newcommand{\ene}{\end{eqnarray}}
\newcommand{\xym}{\xymatrix}
\newcommand{\bp}{\begin{pmatrix}}
\newcommand{\ep}{\end{pmatrix}}

\newcommand{\Hom}{\operatorname{Hom}}

\newcommand{\Aut}{\operatorname{Aut}}

\newcommand{\Id}{\operatorname{Id}}

\newcommand{\PGL}{\operatorname{PGL}}

\def\arraystretch{1.2} 

\begin{document}

\author{W.~D.~Gillam}
\author{A.~Karan}

\thanks{The first author was supported by a Marie Curie/T\"{U}BITAK Co-Funded Brain Circulation Scheme fellowship.}
\address{Department of Mathematics, Bogazici University, Bebek, Istanbul 34342}


\date{October 9, 2015}
\title[The Hausdorff topology as a moduli space]{The Hausdorff topology as a moduli space}

\begin{abstract}  In 1914, F.~Hausdorff defined a metric on the set of closed subsets of a metric space $X$.  This metric induces a topology on the set $H$ of \emph{compact} subsets of $X$, called the \emph{Hausdorff topology}.  We show that the topological space $H$ represents the functor on the category of sequential topological spaces taking $T$ to the set of closed subspaces $Z$ of $T \times X$ for which the projection $\pi_1 : Z \to T$ is open and proper.  In particular, the Hausdorff topology on $H$ depends on the metric space $X$ only through the underlying topological space of $X$.  The Hausdorff space $H$ provides an analog of the Hilbert scheme in topology.  As an example application, we explore a certain quotient construction, called the \emph{Hausdorff quotient}, which is the analog of the Hilbert quotient in algebraic geometry.   \end{abstract}

\maketitle

\arraycolsep=2pt 
\def\arraystretch{1.2}

\section{Introduction}  \label{section:introduction} To fix ideas, consider a compact, metrizable topological space $X$ and a complex projective variety $X'$.  For example, $X$ might be the analytic topological space $(X')^{\rm an}$ underlying $X'$ (cf.\ \cite{Ser}).  The goal of this paper is to construct a topological analog $H$ (the \emph{Hausdorff space} of $X$) of the Hilbert scheme $H'$ of $X'$.  

It will be helpful if we begin with a few words about $H'$ so that our analogy will be more clear.   Throughout the paper, we let $\Sch$ denote the category of schemes of locally finite type over $\CC$ and we refer to objects of $\Sch$ simply as \emph{schemes}.  The scheme $H'$ represents the contravariant functor from $\Sch$ to $\Sets$ taking a scheme $T$ to the set of closed subschemes $Z$ of $T \times X'$ for which the projection $\pi_1 : Z \to T$ is flat.  The set of (closed) points of $H'$ is therefore in bijective correspondence with the set of closed subschemes of $X'$, so that one may view $H'$ as providing some additional structure to this set, determined uniquely by the ``modular interpretation" of $H'$ (the description of the functor that $H'$ represents in terms of closed subschemes of $X'$).  Our goal, then, is to put additional structure (in this case a topology) on the set $H$ of closed subsets of $X$, so that the resulting topological space $H$ has a nice ``modular interpretation" like the one for the Hilbert scheme.  

In fact our $H$ will represent a functor defined in almost the same way as the one represented by $H'$---we just replace ``scheme" with ``space" and ``flat" with ``open".  (The first author has long been a proponent of the idea that open maps of topological spaces are analogous to flat maps of schemes.)  It turns out that the topology on $H$ that will work for our purposes was defined by F.~Hausdorff in his famous book \cite{Hau}.  So the new result here is the modular interpretation of Hausdorff's $H$.  The situation is slightly complicated by the fact that we see no way to define the topology on $H$ ``directly" without choosing a metric for $X$.

Although the book \cite{Hau} appeared over a hundred years ago, the results in this paper could certainly have been obtained by Hausdorff himself at that time.  What stopped that from happening has to do with the history of mathematics.  Before Grothendieck's influence, it would not have been common practice to look for a ``modular interpretation" of such constructions.  In fact, one could imagine an utterly plausible alternative history of algebraic geometry in which the Hilbert scheme $H'$ was constructed long before anyone thought about it as a ``moduli space".  (After all, its construction is largely based on classical commutative algebra of graded rings.)   The fact that history did not play out like that is only because Grothendieck had not only the prescience to think of the right moduli problem (``flat families of closed subschemes") but also the technical expertise to establish its representability.

Like the Hilbert scheme $H'$ of $X'$, the Hausdorff space $H$ of $X$ is, in some sense, a rather large and mysterious object whose main purpose is often to serve as a building block for other constructions.  For example, when an algebraic group $G'$ acts on $X'$, one can use $H'$ to construct a certain quotient of $X'$ by $G'$, called the \emph{Hilbert quotient}.  It is closely related to the similarly-constructed \emph{Chow quotient} and to the quotients considered in geometric invariant theory.  When a topological group $G$ acts on $X$, we can construct an analogous \emph{Hausdorff quotient} by making use of $H$.  The Hausdorff quotient has several advantages over the usual categorical quotient $X/G$---for instance it is always compact Hausdorff.  When $X$ is the underlying analytic space of $X'$, there are some interesting relationships between the Hilbert and Hausdorff quotients.  They often ``agree".  Some instances of the Hausdorff quotient have been examined before---for example, in Morse Theory---though the general construction seems not to have been considered.

\section{The Hausdorff space}  In this section we will give precise statements of our results concerning the Hausdorff space.  Proofs are given in \S\ref{section:proofs}.

Let $X=(X,d)$ be a metric space.  \emph{We assume throughout the paper that for every metric $d$, the distances $d(x,y)$ are bounded above by $1$.}  Since we will be interested only in the topology on $X$ determined by $d$, this assumption is harmless since we can ``cut off" any metric $d$ by replacing $d$ with $\min(d,1)$ without chaning the induced topology on $X$.  Let $C=C(X)$ denote the set of closed subsets of $X$.  For non-empty $A,B \in C$, the \emph{Hausdorff distance} is defined by \be d_C(A,B) & \defeq & \max ( \sup_{a \in A} \inf_{b \in B} d(a,b) , \sup_{b \in B} \inf_{a \in A} d(a,b) ). \ee  Hausdorff showed \cite[Page 293]{Hau} that $d_C$ defines a metric (which is clearly bounded above by $1$) on the set of non-empty $A \in C$.  By declaring the Hausdorff distance between a non-empty closed subset and $\emptyset$ to be $1$ we extend this metric to a metric on $C$.

\begin{example} \label{example:1} Even if two metrics $d$, $d'$ on $X$ determine the same topology on $X$ (hence the same meaning of $C=C(X)$), the corresponding Hausdorff distances $d_C$, $d'_C$ may determine different topologies on $C$.  For example, take $X = \RR^2$, $d$ the usual (but cut off by $1$, as always) metric, and $d'$ the metric on $X$ determined by pulling back the usual (cut off) metric on the open unit disk $D \subset \RR^2$ under the homeomorphism $\RR^2 \to D$ defined by $(x,y) \mapsto (x,y)/(1+|(x,y)|)$.  Let $L_n \subseteq X$ be the line through the origin with slope $n$.  Then $d_C(L_i,L_j)=1$ when $i \neq j$, whereas $L_1, L_2, \dots$ is a Cauchy sequence with respect to $d'_C$ converging to the vertical line through the origin. \end{example}

We shall be primarily interested in the set $H=H(X)$ of \emph{compact} subspaces of $X$.  Of course $H(X)=C(X)$ when $X$ is compact, as was the situation in \S\ref{section:introduction}.  Let $d_H$ be the restriction of $d_C$ to $H$.  The topology on $H$ determined by $d_H$ is called the \emph{Hausdorff topology} and the corresponding topological space $H$ is called the \emph{Hausdorff space} of $X$.  

\begin{thm} \label{thm:A} The subset \be \Z & \defeq \{ (Z,x) \in H \times X : x \in Z \} \ee of $H \times X$ is closed and the projection $\pi_1 : \Z \to H$ is open and proper. \end{thm}

Recall that a subset $Z$ of a topological space $X$ is called \emph{sequentially closed} iff, whenever $z_1,z_2,\dots$ is a sequence of points of $Z$ converging (in $X$) to a point $x$, we have $x \in Z$.  Clearly a closed subset is sequentially closed.  A topological space $T$ is called \emph{sequential} iff every sequentially closed subset of $T$ is closed.  Any first countable space, hence any metrizable space, is sequential.  Let $\SeqTop$ denote the full subcategory of the category $\Top$ of topological spaces consisting of sequential topological spaces.  Our main result is:

\begin{thm} \label{thm:B} The Hausdorff space $H$ of $X$ represents the functor $\SeqTop^{\rm op} \to \Sets$ taking $T \in \SeqTop$ to the set of closed subspaces $Z \subseteq T \times X$ for which the projection $\pi_1 : Z \to T$ is open and proper. \end{thm}

Since the functor in Theorem~\ref{thm:B} clearly depends only on the \emph{topology} on $X$ determined by $d$, the same must be true of the space $H$ representing it, so we have:

\begin{cor} \label{cor:C}  The topology on $H$ determined by the Hausdorff metric $d_H$ depends on the metric $d$ on $X$ only through the topology on $X$ determined by $d$.  In other words, the Hausdorff space $H(X)$ is intrinsically attached to any metrizable space $X$. \end{cor}

(It is not particularly difficult to prove Corollary~\ref{cor:C} directly without using Theorem~\ref{thm:B}.)  Although we do not see how to define the Hausdorff topology on $H$ without choosing a metric, we can still characterize it as follows:

\begin{cor} \label{cor:C2} For a metrizable space $X$, the Hausdorff topology on the set $H$ of compact subsets of $X$ is the smallest topology $\tau$ on $H$ such that the subset $\Z \subseteq H_\tau \times X$ defined in Theorem~\ref{thm:A} is closed and the projection $\pi_1 : \Z \to H_\tau$ is open and proper. \end{cor}

Indeed, if $\tau$ is a topology with those properties, then Theorem~\ref{thm:B} ensures that the identity function defines a continuous map from $H_\tau$ to $H$ with the Hausdorff topology.  It is not obvious that there \emph{exists} a smallest such topology, nor even that there exists \emph{any} such topology.

The following is ``well-known" \cite{Hen} and not due to us in any way, but we mention it for the sake of completeness.

\begin{thm} \label{thm:C} The Hausdorff space of a compact metric space is compact. \end{thm}

Theorem~\ref{thm:C} shows that ``and proper" can be deleted from Theorem~\ref{thm:B} when $X$ is assumed compact.  (This is also easy to see from our proof of Theorem~\ref{thm:B}.)

The Hausdorff space is related to the Hilbert scheme by more than just analogy.  Let us return to the situation mentioned in \S\ref{section:introduction} where $X'$ is a complex projective variety (or, more generally, any projective scheme) and $X$ is its underlying analytic topological space: $X = (X')^{\rm an}$.  Since the ``underlying analytic topological space" functor \be \Sch & \to & \SeqTop \\ X' & \mapsto & (X')^{\rm an} \ee takes closed embeddings to closed embeddings, proper maps to proper maps, and flat maps to open maps \cite{Ser}, the modular interpretations of $H'$ and $H$ yield a continuous map of topological spaces \bne{HilbertHausdorffmap} (H'(X'))^{\rm an} & \to & H((X')^{\rm an}), \ene called the \emph{Hilbert-Hausdorff map}.  It is given by taking a closed subscheme $Z'$ of $X'$ to the underlying analytic closed subspace $(Z')^{\rm an}$ of $(X')^{\rm an}$.  

The map \eqref{HilbertHausdorffmap} is certainly not one-to-one, since a given reduced closed subscheme $Z'$ of $X'$ will typically underlie infinitely many different non-reduced closed subschemes (all of which may even have the same Hilbert polynomial).  The Hausdorff space and the Hilbert scheme also treat ``collisions" differently:  Suppose $Z_1, Z_2, \dots$ (resp. $W_1, W_2, \dots$) is a sequence of reduced closed subschemes of $X'$, each consisting of precisely two (resp. three) points.  In $H'$, the $Z_n$ and the $W_n$ are contained in disjoint open subsets because the $Z_n$ have Hilbert polynomial $2$, whereas the $W_n$ have Hilbert polynomial $3$.  But suppose the $Z_n$ (resp. $W_n$) converge in $H'$ to a (necessarily non-reduced) closed subscheme $Z$ (resp. $W$) supported at a single point $x'$ (resp.\ \emph{the same} $x'$).  Then the closed subspace $\{ x' \}$ of $(X')^{\rm an}$ will be the limit of both the $Z_n^{\rm an}$ and the $W_n^{\rm an}$ in $H$.  Hence the $Z_n^{\rm an}$ and the $W_n^{\rm an}$ won't be contained in disjoint open subsets of $H$, so \eqref{HilbertHausdorffmap} won't be an \emph{embedding} even if we restrict to the subspace of $(H')^{\rm an}$ whose points are the \emph{reduced} closed subschemes of $X'$.

It turns out, however, that the issues discussed above are basically the only thing stopping \eqref{HilbertHausdorffmap} from being an embedding:

\begin{thm} \label{thm:D} Let $X'$ be a projective scheme.  Fix an ample line bundle on $X'$ and a numerical polynomial $p \in \QQ[T]$.  Let $Y$ denote the subspace of $(H'(X'))^{\rm an}$ whose points are \emph{reduced} closed subschemes of $X'$ with Hilbert polynomial $p$ (w.r.t.\ the chosen ample bundle).  Then the restriction of the Hilbert-Hausdorff morphism yields an embedding $Y \into H=H((X')^{\rm an})$ into the Hausdorff space of the underlying analytic space of $X'$. \end{thm}

\section{The Hausdorff quotient}  In this section we will describe a topological analog (the \emph{Hausdorff quotient}) of the \emph{Hilbert quotient} construction in algebraic geometry \cite{BBS}, \cite{Hu}, \cite{Kap1}, \cite[\S1]{Tha}.  We begin by reviewing the latter.  Suppose a group scheme $G'$ acts (algebraically) on a complex projective variety $X'$.  \emph{We assume throughout that our group schemes are finite type over $\CC$ and, for simplicity, connected} (hence irreducible since these group schemes are smooth).  Let $G'$ act on $X' \times X'$ by acting on the second factor.  The $G'$ orbit $V$ of the diagonal $X' \cong \Delta \subseteq X' \times X'$ is an irreducible, constructible subset of $X' \times X'$ because it is the image of the irreducible space $G' \times \Delta$ under the action morphism.  We regard $V$ and its closure $Z$ in $X' \times X'$ as integral schemes of finite type over $\CC$ by giving them the reduced-induced scheme structures from $X' \times X'$.  By using ``generic flatness" results (cf.\ \cite[6.9]{EGA}) and basic properties of flat, proper maps (cf.\ \cite[12.2.4]{EGA}) we can find a non-empty open subspace $U'$ of $X'$ which is \emph{stable} in the sense that $Z \cap (U' \times X')$ is flat over $U'$ (via $\pi_1$) and the fiber $Z_{x'}$ of $Z$ over any point $x' \in U'$ is the orbit closure $\ov{G' x'}$ (equivalently: the open subset $G'x'=V_{x'} \subseteq Z_{x'}$ is dense).  From the modular interpretation of the Hilbert scheme $H'$ of $X'$ we obtain a morphism $e'_{U'} : U' \to H'$, given on closed points by taking $x' \in U'$ to $\ov{G'x'}$.  The Hilbert quotient $X'/_{H'} G'$ is defined to be the closure of the image of $e'_{U'}$ (this is clearly independent of the choice of stable $U'$), with its reduced-induced scheme structure from $H'$.  It is projective.

Now let's look at the topological analog where $G$ is a topological group acting continuously on a metrizable space $X$.  Write $H=H(X)$ for the Hausdorff space of $X$.  We assume the orbit closure $\ov{Gx}$ is compact for every $x \in X$.  (This holds automatically if $G$ or $X$ is compact.)  We write $\UU$ for the set of open, dense, $G$-invariant subsets of $X$.  We have a map \emph{of sets} \be e : X & \to & H \\ x & \mapsto & \ov{Gx}, \ee which is \emph{not generally a continuous map of topological spaces}.  Never-the-less, for any $U \in \UU$, we can form the closure $\ov{e(U)}$ of the image $e(U)$ of $U$ in $H$.  The \emph{Hausdorff quotient} of $X$ by $G$, denoted $X /_H G$, is defined by \bne{Hausdorffquotient} X /_H G & \defeq & \bigcap_{U \in \UU} \ov{e(U)}. \ene  By construction $X /_H G$ is a closed subspace of $H$, so it is metrizable.  It is also compact when $X$ is compact (by Theorem~\ref{thm:C}).

\begin{rem} \label{rem:Hausdorffquotient1}  The construction of the Hausdorff quotient $X /_H G$ does not make any use of the topology of $G$ or the continuity of the action, except possibly in so far as compactness of $G$ and continuity of the action would be one way to ensure compactness of the orbit closures.  However, assumptions of this nature become important if one wants to prove anything interesting about $X /_H G$. \end{rem}

\begin{rem} \label{rem:Hausdorffquotient2} In fact, the compactness of the orbit closures is not really necessary to define the Hausdorff quotient, since one could use $C(X)$ instead of $H(X)$.  But then the dependence of $X /_H G$ on the choice of metric inducing the topology of $X$ becomes an issue (cf.\ Example~\ref{example:1}). \end{rem}

\begin{defn} \label{defn:stable} A set $U \in \UU$ is called \emph{stable} iff $e|U : U \to H$ is continuous.  A set $U \in \UU$ is called \emph{semi-stable} iff the natural inclusion $X /_H G \subseteq \ov{e(U)}$ of closed subspaces of $H$ is an equality. \end{defn}

It is clear that stable implies semi-stable and that stability (resp.\ semi-stability) of $U \in \UU$ implies stability (resp.\ semi-stability) of any $V \in \UU$ with $V \subseteq U$.

\begin{example} The group $G=\RR_{>0}$ of positive real numbers under multiplication acts continuously on the closed interval $X=[0,\infty]$ by multiplication.  We have \be \UU & = & \{ [0,\infty], (0,\infty], [0,\infty), (0,\infty) \}. \ee  In this case ``stable" and ``semi-stable" are equivalent and only $(0,\infty) \in \UU$ has this property.  The Hausdorff quotient $X/_H G$ is a point. \end{example}

\begin{thm} \label{thm:E} If $G$ and $X$ are compact, then the Hausdorff quotient $X /_H G$ coincides with the usual topological quotient $X/G$. \end{thm}

\begin{rem} Theorem~\ref{thm:E} probably holds under weaker assumptions on $X$ and $G$. \end{rem}

\begin{example} Let $X$ be a compact Riemannian manifold without boundary, $f : X \to \RR$ a Bott-Morse function on $X$.  We get an action of $G=(\RR,+)$ on $X$ by integrating the gradient vector field of $f$.  The Hausdorff quotient $X/_H G$ ``should be" (one might say ``defines") the moduli space of (possibly broken) flow lines of $f$ (from the ``top" critical locus to the ``bottom" critical locus).  See \cite{Gil} for further discussion. \end{example}

Now suppose the topological situation ($G$ acting on $X$) is the underlying analytic topological picture of the algebraic situation ($G'$ acting on $X'$).  Then for any stable (in the algebraic sense) $U' \subseteq X'$, the map (a priori \emph{of sets}) $e|(U')^{\rm an} : (U')^{\rm an} \to H$ is the composition of the continuous map $(e'_{U'})^{\rm an} : (U')^{\rm an} \to (H')^{\rm an}$ and the Hilbert-Hausdorff map $(H')^{\rm an} \to H$, so it \emph{is} continuous, hence $(U')^{\rm an}$ is stable (in the topological sense).  It follows that the Hilbert-Hausdorff map $(H')^{\rm an} \to H$ takes the analytic topological space $(X'/_{H'} G')^{\rm an} \subseteq (H')^{\rm an}$ underlying the Hilbert quotient into the Hausdorff quotient $X/_H G \subseteq H$.  We thus obtain a \emph{Hilbert-Hausdorff quotient map} \bne{HHquotientmap} (X'/_{H'} G')^{\rm an} & \to & X/_H G. \ene  Using Theorem~\ref{thm:D}, one easily proves:

\begin{thm} \label{thm:F} Let $X'$ be a projective variety equipped with an action of a group scheme $G'$.  Assume every (closed) point $Z'$ of the Hilbert quotient $X'/_{H'}G' \subseteq H'$ is a \emph{reduced} closed subscheme of $X'$.  Then the Hilbert-Hausdorff quotient map \eqref{HHquotientmap} is a homeomorphism. \end{thm}

\begin{example} The Hilbert quotient of $(\PP^1)^n$ by the diagonal action of $\PGL_2 = \Aut( \PP^1 )$ is the moduli space $\ov{\M}_{0,n}$ of stable, marked genus zero curves (see \cite{Kap1}, \cite{Kap2}, \cite{GG}).  The hypothesis of Theorem~\ref{thm:E} is satisfied---the reduced closed subscheme of $(\PP^1)^n$ corresponding to a point of the Hilbert quotient $\ov{\M}_{0,n}$ is described explicitly in \cite{GG}.  It follows that the Hausdorff quotient of the topological space $(\CC \PP^1)^n$ by the diagonal action of $\PGL_2(\CC)$ is the underlying analytic topological space of $\ov{\M}_{0,n}$. \end{example}

\section{Proofs} \label{section:proofs}  We denote the open ball of radius $\epsilon > 0$ centered at a point $x$ of a metric space $X=(X,d)$ by $B(x,\epsilon)$.  We always use lower case letters for points of $X$ and capital letters for points of $H=H(X)$, so an expression like $B(x,\epsilon)$ (resp.\ $B(Z,\epsilon)$) clearly refers to a ball in $(X,d)$ (resp.\ $(H,d_H)$).  The formula \be \rho((Z,x),(Z',x')) & \defeq & \max( d_H(Z,Z'), d(x,x') ) \ee defines a metric $\rho$ on $H \times X$ inducing the product topology.  We always use the metric $\rho$ when refering to a ball $B((Z,x),\epsilon)$ in $H \times X$.  

\subsection{Proof of Theorem~\ref{thm:A}}  To show that $\Z$ is closed in $H \times X$, suppose $(Z,x) \notin \Z$.  Then $x \notin Z$ and $Z \subseteq X$ is closed, so there is an $\epsilon >0$ such that $B(x,\epsilon) \cap Z = \emptyset$.  We claim that $B((Z,x),\epsilon/2) \cap \Z = \emptyset$.  If not, then there is $(Z',x') \in B((Z,x),\epsilon/2)$ with $x' \in Z'$.  By definition of $B((Z,x),\epsilon/2)$ we have $d_H(Z,Z') < \epsilon/2$ and $d(x,x') < \epsilon/2$.  By definition of $d_H$, the first inequality implies $d(z,x') < \epsilon/2$ for some $z \in Z$.  The triangle inequality then gives $d(x,z) < \epsilon$, contradicting $B(x,\epsilon) \cap Z = \emptyset$.  (Note that we do not use compactness of $Z$ anywhere, so this same argument shows that $\{ (Z,x) \in C(X) \times X : x \in Z \}$ is closed in $C(X) \times Z$.)

To show that $\pi_1 : \Z \to H$ is proper, first note that the fiber of $\pi_1$ over $Z \in H$ ``is" $Z$, which is compact by definition of $H$, so $\pi_1$ is a map between Hausdorff spaces with compact fibers, so to show that it is proper it remains only to show that it is closed.  This boils down to the following:  Whenever $Z_1, Z_2, \dots$ is a sequence in $H$ converging to $Z \in H$ and we have points $x_1 \in Z_1, x_2 \in x_2, \dots$, then there is a point $x \in Z$ in the closure of $\{ x_1, x_2, \dots \}$ in $X$.  Since the $Z_i$ converge to $Z$, we can assume, after possibly passing to a subsequence, that there are points $y_1,y_2, \dots \in Z$ with $d(x_n,y_n) < 1/n$.  Since $Z$ is compact, we can assume, after possibly passing to a subsequence, that $y_1, y_2, \dots $ converges to $x \in Z$.  To see that this $x$ is in the closure of $\{ x_1, x_2, \dots \}$, fix $\epsilon > 0$.  Choose $n$ large enough that $d(x_n,y_n) < 1/n < \epsilon/2$ and $d(y_n,x) < \epsilon / 2$.  Then $x_n \in B(x,\epsilon)$ by the triangle inequality.

To show that $\pi_1 : \Z \to H$ is open, it suffices to show that \be \pi_1(\Z \cap B((Z,x),\epsilon)) & = & B(Z,\epsilon) \ee for any $(Z,x) \in \Z$ and any $\epsilon > 0$.  The containment $\subseteq$ is obvious.  For the opposite containment, consider a point $Z' \in B(Z,\epsilon)$.  By definition of $d_H$ the inequality $d_H(Z,Z')< \epsilon$ ensures that $Z \subseteq \cup_{x' \in Z'} B(z',\epsilon)$.  Since $x \in Z$, there is hence some $x' \in Z'$ such that $x \in B(x',\epsilon)$.  Hence $x' \in B(x,\epsilon)$ and therefore $Z' = \pi_1(Z',x') \in \pi_1( \Z \cap B((Z,x),\epsilon))$, as desired.  (Again note that we do not use compactness anywhere, so this same argument shows that the projection from $\{ (Z,x) \in C(X) \times X : x \in Z \}$ to $C(X)$ is open.)

\subsection{Proof of Theorem~\ref{thm:B}}  We first prove:

\begin{thm} \label{thm:Bstrong} Let $T$ be a topological space.  \begin{enumerate} \item \label{Bstrong1} If $f : T \to H$ is continuous, then $$Z_f \defeq \Z \times_T H = \{ (t,x) \in T \times X : (f(t),x) \in \Z \} = \bigcup_{t \in T} \{ t \} \times f(t) $$ is a closed subspace of $T \times X$ and the projection $\pi_1 : Z_f \to T$ is open and proper. \item \label{Bstrong2} If $T$ is sequential and $Z$ is a closed subspace of $T \times X$ for which the projection $\pi_1 : Z \to T$ is open and proper then $Z_t \defeq \{ x \in X: (t,x) \in Z \}$ is a compact subset of $X$ for each $t \in T$ and the function $f_Z : T \to H$ defined by $f_Z(t) \defeq Z_t$ is continuous. \end{enumerate} \end{thm}

Assuming Theorem~\ref{thm:Bstrong}, we see that for $T \in \SeqTop$, the map $f \mapsto Z_f$ establishes a bijection (natural in $T$) from $\Hom_{\Top}(T,H)$ to $\F(T)$ ($\F$ the functor described in Theorem~\ref{thm:B}) with inverse $Z \mapsto f_Z$, thus proving Theorem~\ref{thm:B}.

Theorem~\ref{thm:Bstrong}\eqref{Bstrong1} follows from Theorem~\ref{thm:A} because closed embeddings, open maps, and proper maps are all stable under base change and, when $f$ is continuous, \bne{cartdiagram} \xym{ Z_f \ar[r]^-{\subseteq} \ar[d] & T \times X \ar[d]^{f \times \Id} \ar[r]^-{\pi_1} & T \ar[d]^f \\ \Z \ar[r]^-{\subseteq} & H \times X \ar[r]^-{\pi_1} & H } \ene is a cartesian diagram of topological spaces.

To prove Theorem~\ref{thm:Bstrong}\eqref{Bstrong2}, we will make use of the following lemma, whose proof is left to the reader as an exercise with the definition of ``sequential space."

\begin{lem} \label{lem:thelemma} Let $\ov{\NN} = \NN \cup \{ \infty \}$ denote the one point compactification of the countably infinite discrete space $\NN = \{ 1,2, \dots \}$.  Let $T$ be a sequential space, $H$ an arbitrary topological space, $f : T \to H$ an arbitrary map of sets.  Then $f$ is continuous iff the composition $fg$ is continuous for every continuous map $g : \ov{\NN} \to T$. \end{lem}

Since the hypotheses on $Z \subseteq T \times H$ and $\pi_1 : Z \to T$ in Theorem~\ref{thm:Bstrong}\eqref{Bstrong2} are stable under base change along any map $g : T' \to T$ and the map $f_Z g : T' \to H$ agrees with the map $f_{Z'}$ constructed from the pullback $Z' \defeq Z \times_T T'$, Lemma~\ref{lem:thelemma} reduces Theorem~\ref{thm:Bstrong}\eqref{Bstrong2} to the case where $T = \ov{\NN}$.  In this case, continuity of $f_Z : \ov{\NN} \to H$ is equivalent to continuity at $\infty$, which, from the definition of $d_H$, is easily seen to be equivalent to:  For every $\epsilon > 0$, there is a neighborhood $U$ of $\infty$ in $\ov{\NN}$ satisfying: \begin{enumerate} \item \label{one} $Z_n \subseteq \cup_{z \in Z_\infty} B(z,\epsilon)$ for all $n \in U$ and \item \label{two} $Z_\infty \subseteq \cup_{z \in Z_n} B(z,\epsilon)$ for all $n \in U$. \end{enumerate}   It is enough to find neighborhoods $U_1$ (resp.\ $U_2$) (of $\infty$ in $\ov{\NN}$) satisfying \eqref{one} (resp.\ \eqref{two}), for then we can take $U = U_1 \cap U_2$.

Suppose no such $U_1$ exists.  Then after possibly replacing $\NN=\{ 1,2,\dots \}$ with some subsequence, we can find points $x_n \in Z_n$ (i.e.\ points $(n,x_n) \in Z \subseteq \ov{\NN} \times X$) such that \bne{ineq} d(x_n,Z_\infty) & \geq & \epsilon \quad \quad \forall n \in \NN. \ene  Since $\ov{\NN}$ is compact Hausdorff and $\pi_1 : Z \to \ov{\NN}$ is proper (by assumption), $Z$ is also compact Hausdorff, so we can assume, after passing to a subsequence, that the points $(n,x_n)$ of $Z$ converge to a point $(\infty,x) \in Z$.  Then $x \in Z_\infty$ and the $x_n$ converge to $x$ in $X$, contradicting \eqref{ineq}.

To find $U_2$ as desired, we use compactness of $Z_\infty$ (which holds since $Z_\infty$ is a fiber of the proper map $\pi_1 : Z \to T$) to find $x_1, \dots, x_n \in Z_\infty \subseteq X$ such that \bne{containment} Z_\infty & \subseteq & \cup_{i=1}^n B(x_i,\epsilon/2). \ene  We claim that \be U_2 & \defeq & \cap_{i=1}^n \pi_1( Z \cap (T \times B(x_i,\epsilon/2))) \ee (which is open by the assumption that $\pi_1 : Z \to T$ is open) is as desired.  To see this, fix some $n \in U_2$ and an arbitrary point $x \in Z_\infty$.  By \eqref{containment}, there is some $i \in \{ 1, \dots, n \}$ such that $d(x,x_i) < \epsilon/2$.  Since $U_2 \subseteq \pi_1( Z \cap (T \times B(x_i,\epsilon/2)))$, there is a $y \in Z_n$ with $d(x_i,y) < \epsilon / 2$.  Then $d(x,y)<\epsilon$ by the triangle inequality.  This proves that $Z_\infty \subseteq B(Z_n,\epsilon)$, as desired.

\subsection{Proof of Theorem~\ref{thm:D}}  Let $H'_p$ be the open and closed subscheme of the Hilbert scheme $H'$ of $X'$ parametrizing closed subschemes of $X'$ with Hilbert polynomial $p$ (w.r.t.\ the chosen ample bundle), so that $Y$ is a subspace of $(H'_p)^{\rm an}$.  Denote the restriction of the Hilbert-Hausdorff morphism \eqref{HilbertHausdorffmap} to $(H'_p)^{\rm an}$ by $f : (H'_p)^{\rm an} \to H$.  Since $H'_p$ is a projective (hence proper) scheme, $(H'_p)^{\rm an}$ is a compact Hausdorff space.  Since $H$ is Hausdorff, $f$ is a closed (and continuous) map, so to prove that its restriction $g \defeq f|Y : Y \to H$ is an embedding it will be enough to check that $g$ is one-to-one and that \bne{desiredequality} W \cap Y & = & g^{-1}(f(W)) \ene for any closed subset $W$ of $(H'_p)^{\rm an}$.  The map $g$ is one-to-one because two reduced closed subschemes of $X'$ with the same underlying analytic space in particular have the same set of closed points, hence also the same underlying Zariski topological space, hence are equal.  The containment $\subseteq$ in \eqref{desiredequality} is obvious and the only issue for the reverse containment is to show that if $g(Z) = f(Z')$ for some $Z \in Y$, $Z' \in (H'_p)^{\rm an}$, then $Z=Z'$.  The point here is that since $Z$ is reduced and the closed subscheme $Z'$ has the same underlying analytic (hence Zariksi) space as $Z$, we have a scheme-theoretic containment $Z \subseteq Z'$.  But since $Z$ and $Z'$ have the same Hilbert polynomial $p$, this containment must be an equality.

\subsection{Proof of Theorem~\ref{thm:E}} Since $X$ is compact Hausdorff and $G$ is compact, the quotient map $q : X \to X/G$ is proper (one does not even need $G$ to be Hausdorff for this) and every $G$ orbit $Gx$ is compact Hausdorff.  Like the quotient map for any group action, $q$ is also open and surjective.  The map $e : X \to H$ defined by $e(x) \defeq Gx$ is hence continuous by Lemma~\ref{lem:gen} below.  That is, $X$ itself is stable, so that $X/_H G$ is the closure of $e(X)$ in $H$.  Since $e$ is clearly constant on $G$ orbits it factors (by the universal property of $q$) through $q$ via a continuous map $\ov{e} : X/G \to H$.  Lemma~\ref{lem:gen} shows that $\ov{e}$ is an embedding (hence a closed embedding because its domain $X/G$ is compact).  Since $\ov{e}(X/G)=e(X)$, the proof is complete.

\begin{lem} \label{lem:gen} Let $X$ be a metrizable space with Hausdorff space $H$, $q : X \to Y$ an open, proper, surjective map.  Then the function $e : Y \to H$ defined by $e(y) \defeq q^{-1}(y)$ is a continuous embedding. \end{lem}

To prove this, first note that $Y$ is sequential:  If $Z \subseteq Y$ is sequentially closed, then $q^{-1}(Z)$ is sequentially closed in $X$ since $q$ is continuous, hence closed since $X$ is metrizable, hence $Z=q(q^{-1}(Z))$ is closed since $q$ is closed and surjective.  Next note that the ``graph" map $x \mapsto (f(x),x)$ defines a homeomorphism from $X$ to a closed subspace $Z$ of $Y \times X$ identifying $q$ with $\pi_1 : Z \to Y$, so that $e$ is nothing but the continuous map obtained from $Z$ and the universal property of the Hausdorff quotient.  Clearly $e$ is one-to-one since $q$ is surjective.  To see that it is an embedding, one checks that for any open subset $U \subseteq Y$, the subset \be \tilde{U} & \defeq & \{ W \in H : W \subseteq q^{-1}(U) \} \ee of $H$ is open and satisfies $U=e^{-1}(\tilde{U})$.

\end{document}